\documentclass[12pt]{amsart}
\usepackage{amssymb,amsthm,amsmath,a4}

\newcommand{\area}{{\rm Area}}
\newtheorem{theorem}[equation]{Theorem}

\theoremstyle{remark}
\newtheorem{remark}[equation]{Remark}
\theoremstyle{definition}

\numberwithin{equation}{section}
\begin{document}
\title[Pleijel's nodal domain theorem]{Pleijel's nodal domain theorem \\ for free membranes}
\author{Iosif Polterovich}
\address{D\'e\-par\-te\-ment de math\'ematiques et de
sta\-tistique, Univer\-sit\'e de Mont\-r\'eal, CP 6128 succ.
Centre-Ville, Mont\-r\'eal,  QC  H3C 3J7, Canada.}
\email{iossif@dms.umontreal.ca}
\thanks{Research partially supported by NSERC and FQRNT}
\begin{abstract}
We prove an analogue of Pleijel's nodal domain theorem for piecewise
analytic planar domains with Neumann boundary conditions. This
confirms a conjecture made by Pleijel in 1956. The proof is a
combination of Pleijel's original approach and an estimate due to
Toth and Zelditch for the number of boundary zeros of Neumann
eigenfunctions.
\end{abstract}
\maketitle
\section{Introduction}
Let $\Omega \subset \mathbb{R}^2$ be a bounded planar domain. Let
$0\le \lambda_1 < \lambda_2\le \dots \le \lambda_k\le \dots$ be the
eigenvalues of the Laplacian on $\Omega$ with either Dirichlet or
Neumann boundary conditions, and $\phi_1,\phi_2,\dots, \phi_k,\dots$
be an orthonormal basis of eigenfunctions: $\Delta \phi_k=\lambda_k
\phi_k$, $k=1,2,\dots$. We assume that the boundary $\partial
\Omega$ is sufficiently regular (say, piecewise smooth). In the
Neumann case we also assume throughout the paper that the internal
cone condition \cite[chapter V.4]{EE} is satisfied at the corners in
order to ensure that the spectrum is discrete. According to the
classical Courant's nodal domain theorem \cite{Courant} (see also
\cite[p. 452]{CH}), the number $n_k$ of nodal domains of the %$k$-th
eigenfunction $\phi_k$ %of the Laplacian
is at most $k$. In 1956, Pleijel \cite{Pl} showed that for planar
domains with Dirichlet boundary conditions Courant's bound can be
asymptotically improved:
\begin{equation}
\label{bound:pleijel} \limsup_{k \to \infty} \frac{n_k}{k} \le
\frac{4}{j_1^2} \approx 0.691...,
\end{equation}
where $j_1\approx 2.4$ is the first zero of the Bessel function
$J_0$. In particular, this implies that the equality in Courant's
theorem is attained only for a finite number of eigenfunctions. Note
that in the $1$-dimensional case, the eigenfunction $\phi_k$ has
{\it exactly} $k$ nodal domains for any $k=1,2,\dots$.

The proof of the estimate \eqref{bound:pleijel} is an application of
the Faber--Krahn inequality for the first eigenvalue (\cite{F,K},
see also \cite[p. 87]{Chavel}). In \cite[section 7]{Pl}, Pleijel
writes regarding the inequality \eqref{bound:pleijel} that {\it
``...it seems highly probable that the result... is also true for
free membranes.''} Recall that a free membrane corresponds to the
Neumann boundary value problem. The difficulty in this case is that
one can not apply the Faber--Krahn inequality to nodal domains that
are adjacent to $\partial \Omega$: they have Neumann conditions on a
part of their boundary (here and further on we say that a nodal
domain $D$ is {\it adjacent} to $\partial \Omega$ if ${\rm
length}\,(\partial D\, \cap\,
\partial \Omega)>0$). Later, Pleijel's result was generalized in
\cite{Pe, BM, Berard} to higher--dimensional domains with Dirichlet
boundary conditions and to compact closed manifolds. However, the
case of domains with Neumann boundary conditions remained open.

The purpose of this note is to prove an analogue of Pleijel's
theorem for piecewise analytic planar domains with Neumann boundary
conditions. We consider separately the nodal domains that are
adjacent to $\partial \Omega$, and  the remaining nodal domains
having pure Dirichlet conditions on their boundaries. In the latter
case, we use the original Pleijel's approach, while in the former we
apply an estimate recently obtained by Toth--Zelditch \cite{TZ} for
the number of boundary zeros of Neumann eigenfunctions. In fact, it
follows from \cite[Theorem 2]{TZ} that the number of nodal domains
adjacent to $\partial \Omega$ is small compared to
Courant's bound, %and therefore its contribution to \eqref{bound:pleijel} is negligible,
see~\eqref{bound:boundary}.
\begin{theorem}
\label{main} Let $\Omega \subset \mathbb{R}^2$ be a piecewise real
analytic domain with Neumann boundary conditions. Then the estimate
\eqref{bound:pleijel} holds.
\end{theorem}
\section{Proof of Theorem \ref{main}}
Let $m_k$ be the number of nodal domains of the eigenfunction
$\phi_k$ that are adjacent to $\partial \Omega$ and $l_k=n_k-m_k$ be
the number of the remaining nodal domains. We show first that
\begin{equation}
\label{bound:boundary} \limsup_{k \to \infty} \frac{m_k}{k}=0
\end{equation}
Let $\mathcal{Z}_k=\{x\in \partial \Omega : \phi_k(x)=0 \}$.
According to \cite[Theorem 2]{TZ}, the number $N_k~=~{\rm
card}(\mathcal{Z}_k)$ of boundary zeros of the eigenfunction
$\phi_k$ is bounded by $C\sqrt{\lambda_k}$, where $C>0$ is a
constant depending only on $\Omega$. It is easy to check that
$m_k~\le~N_k+s$, where $s$ is the number of connected components of
$\partial \Omega$. Indeed, the set $\partial \Omega\setminus
\mathcal{Z}_k$ has at most $N_k+s$ connected components, and for
each connected component $\Gamma \subset
\partial \Omega\setminus \mathcal{Z}_k$ there exists a unique nodal
domain $D$ adjacent to the boundary such that $\Gamma \subset
\partial D$ (note that this correspondence is not necessarily one-to-one:
the boundary of a nodal domain may contain more than one connected
component of $\partial \Omega\setminus \mathcal{Z}_k$). At the same
time, by Weyl's law (see \cite[p.~9]{Chavel})
\begin{equation}
\label{Weyl} \lim_{k \to \infty} \frac{\lambda_k}{k}=\frac{4\pi
}{\area(\Omega)}.
\end{equation}
Therefore, for some constant $C_1>0$,
$$
\limsup_{k\to\infty} \frac{m_k}{k} \le \limsup_{k \to \infty}
\frac{N_k+s}{k} \le \limsup_{k\to \infty}
\frac{C_1\sqrt{\lambda_k}}{\lambda_k}=0,
$$
and this completes the proof of \eqref{bound:boundary}.

Now let us show that
\begin{equation}
\label{bound:interior} \limsup_{k \to \infty} \frac{l_k}{k} \le
\frac{4}{j_1^2}.
\end{equation}
The proof of this bound is similar to the proof of
\eqref{bound:pleijel}, see \cite[pp. 24-25]{Chavel}. Indeed, let
$\tilde \Omega \subset \Omega$ be the union of all nodal domains
that are not adjacent to~$\partial \Omega$. Every such domain $D$
has pure Dirichlet conditions on the boundary, and therefore one can
apply the Faber-Krahn inequality:
$$
\lambda_k(\Omega)\area(D)=\lambda_1(D)\area(D) \ge \pi j_1^2.
$$
Summing up the left and the right hand sides over all domains $D$
that are not adjacent to the boundary we get:
$$
\lambda_k(\Omega) \area (\tilde \Omega) \ge l_k \pi j_1^2.
$$
Therefore, taking into account that $\area(\tilde \Omega) < \area
(\Omega)$ and using  \eqref{Weyl}, one obtains
$$
\limsup_{k\to \infty} \frac{l_k}{k} \le \limsup_{k\to\infty}
\frac{\lambda_k(\Omega) \area(\tilde \Omega)}{\pi j_1^2 k} \le
\frac{4}{j_1^2}.
$$
Since $n_k=m_k+l_k$, putting together \eqref{bound:boundary} and
\eqref{bound:interior} we complete the proof of the theorem. \qed

\begin{remark}
Theorem \ref{main} is proved under the same assumptions as
\cite[Theorem~2]{TZ}. It would be interesting to extend Theorem
\ref{main} to higher dimensions and to replace piecewise analyticity
by a weaker condition. It is likely that the results of \cite{TZ},
and hence Theorem \ref{main} as well, hold also for domains with
mixed Dirichlet--Neumann boundary conditions.
\end{remark}
\begin{remark}
It is clear from the proof of \eqref{bound:pleijel} that this
estimate is {\it not} sharp for both Dirichlet  and Neumann boundary
conditions. Indeed, the Faber-Krahn inequality is an equality only
for the disk, and nodal domains of an eigenfunction can not be {\it
all} disks at the same time. Therefore, a natural question is to
find an optimal constant in \eqref{bound:pleijel}. Motivated by  the
results of \cite{BGS} we suggest that for any regular bounded planar
domain with either Dirichlet or Neumann boundary conditions,
$\limsup_{k \to \infty} \frac{n_k}{k} \le \frac{2}{\pi} \approx
0.636...$. If true, this estimate (which is quite close to Pleijel's
bound) is sharp and attained for the basis of separable
eigenfunctions on a rectangle ~\cite{Pl, SS}.
\end{remark}

\end{document}